\documentclass[11pt,oneside]{amsart}

\usepackage{fullpage}
\usepackage{psfrag}
\usepackage{setspace}
\usepackage{paralist}

\usepackage{color}
\usepackage{array,amscd,amsmath,amsthm,ifthen}
\usepackage{amssymb,pstricks,amscd,epsfig}
\usepackage{amssymb}
\usepackage{graphicx}
\usepackage[all]{xy}
\usepackage{multicol}
\usepackage{cite}
\usepackage{enumitem}

\usepackage{comment}

\numberwithin{equation}{section}

\usepackage[colorlinks,allcolors=red]{hyperref}

\def\I{\mathcal{I}}

\def\cO{\mathcal{O}}

\def\PP{\mathbb{P}}

\def\ov#1{\overline{#1}}

\def\cliff{\mathrm{Cliff}}
\def\Pic{\operatorname{Pic}}

\def\supp{\operatorname{Supp}}

\def\rk{\operatorname{rk}}

\def\cliff{\operatorname{Cliff}}
\def\gon{\operatorname{gon}}
\def\ind{\operatorname{index}}

\newtheorem{thm}{Theorem} [section]
\newtheorem{prop}[thm]{Proposition}
\newtheorem{lemma}[thm]{Lemma}
\newtheorem*{claim*}{Claim}

\newtheorem{rem}[thm]{Remark}
\newtheorem{defin}[thm]{Definition}

\newcommand{\be}{\begin{equation}}
\newcommand{\ee}{\end{equation}}

\begin{document}

\author{Enrico Arbarello}

\email{enrico.arbarello@gmail.com}

\title {A remark  on  du Val  linear systems.}

\begin{abstract} Let $|L_g|$, be the  genus $g$ du Val linear system on a Halphen surface $Y$ of index $k$. We prove that the Clifford index $\cliff(C)$ is constant on smooth
curves $C\in |L_g|$. 
Let $\gamma(C)$ be the gonality of $C$.
 When $\cliff(C)<\lfloor{\frac{g-1}{2}}\rfloor$ (the  relevant case), we show that  $\gamma(C)=\cliff(C)+2=k$, 
and that the gonality is realized by the Weierstrass linear series  $|-{kK_Y}_{|C}|$, which is totally ramified at one point.
The proof of the first statement follows closely the path indicated by Green and Lazarsfeld  for a similar statement regarding K3 surfaces.  \end{abstract}

\maketitle

\section{introduction}

Halphen surfaces and du Val curves came to light  accidentally in \cite{ABS17},  among a long list of 
possible singular limits of polarized K3 surfaces. The similarities of these limits 	 with their nonsingular pre-limit counterparts were 
studied  in 
 \cite{ABFS16}, \cite{AB17}, 
  the terminology  ``{\it du Val linear system}"  was introduced in\cite{ABFS16}. 
In this utilitarian note we push a little bit further the strong analogy between du Val linear systems and the system of hyperplane sections 
of  polarized K3 surfaces, by looking at 
Clifford index and  gonality (see Theorem \ref{main}).
\vskip 0.2 cm
It is a pleasure  to  thank Andrea Bruno, Margherita Lelli-Chiesa, Giulia Sacc\`a, and Edoardo Sernesi for useful exchanges on the subject of this note.

\section{ du Val linear systems}

%%%
%%%
From  \cite{C-D12}, we recall the following definition:
$$
\text{\parbox{.88\textwidth}{{\it A smooth rational projective surface X is a Halphen surface if there exists an integer $k>0$ such that the linear system  $|-kK_Y|$ is one-dimensional, has no fixed component, and  no base point. The index of a Halphen surface is the smallest possible value for such a positive integer $k$. }}}
$$
Halphen surfaces are obtained by blowing up  $\PP^2$ at nine points, and contain no $(-n)$ curves (i.e. smooth rational curves of self-intersection equal to $-n$), for $n\neq1,2$. Moreover, there is a precise set of open  conditions on the set of nine points  
for a Halphen surface, of given index $k$, to be unnodal, i.e.  not containing  $(-2)$-curves, e.g.: the nine points should be distinct, no
three of them on a line, and so on (cf. \cite [Proposition 2.5] {C-D12}).
We fix, once and for all, an unnodal Halphen surface $Y$ of index $k$, obtained as a blow-up 
 $$
 \pi: Y\to \PP^2
 $$ 
 at nine distinct points $p_1,\dots, p_9$. We denote
 by $E_1,\dots, E_9$ the exceptional divisors of the blow-up, and by $\ell$ the class of the
 inverse image of a line in $\PP^2$. We may and  will assume that there is a unique (and smooth) cubic  through the nine points, and we will denote   by $J_Y\subset Y$ the proper transform of this cubic under the blow up, so that

$$
K_Y=-J_Y
$$
Therefore, the index of the Halphen surface $Y$  is the minimal  integer $k$ such that
$h^0(Y, kJ_Y)=2$, and $h^0(Y, (k-1)J_Y)=1$. We have

\be\label{characters}
J_Y=3\ell-E_1-\cdots-E_9\,,\qquad {J_Y}^2=0
\ee
Next, we define the {\it du Val linear system $L_g$ on $Y$} by setting:
\be\label{duvY}
 L_g=|3g\ell-gE_1-\cdots-gE_8-(g-1)E_9|\,,\,\,\text{and also}\,\,\,L_g=|C'|
\ee 
We often denote by $C'$ a smooth element of $L_g$. 
From the point of view of plane curves, the smooth elements of this linear system correspond to degree-$3g$  plane curves having eight points of
multiplicity $g$ and one point of multiplicity $g-1$, and no other singularities. As we shall presently recall, one of the main features of  these plane curves is that they behave like hyperplane sections of an elliptic genus-$g$ polarized K3 surface. The ``elliptic" attribute is due to the ubiquitous presence of the elliptic pencil $E=|-kK_Y|$.
Let us enumerate a few facts, (see for instance \cite[Section 2]{ABFS16}):

{\rm

\begin{enumerate}[label=(\roman*)]

 \item $\dim L_g=g$. This follows from the Riemann-Roch Theorem, since for  $C'\in L_g$ we have $h^1(Y,C')=1$.

 \item If the index of  $Y$ is greater than or equal to 3, the general element of $L_g$ is smooth \cite{Harb85}.

 \item$L_{g-1}=|L_g-J_Y|$. This key property  is at the basis of the induction procedure used in the recent and  quite remarkable paper \cite{B-LC23}.

 \item $C'\cdot J_Y=1\,,\,\, C'\cap J_Y=p_{10}(g)$, so that $p_{10}(g)$  is the unique base point of $L_g$.
The point $p_{10}(g)\in J_Y$ depends only on $Y$ and $g$, and plays an important role.

\item By the Hodge Index Theorem (HIT), given any decomposition of $L_g$ into a sum of effective divisors $|D+D'|=L_g$,  with  
$D\cdot J_Y=0$, then $D=\alpha J_Y$, for some positive $\alpha$ (\cite{Nag60}, and   also \cite[Proposition 2.3]{ABFS16}).

\item If $Y$ has index  $g+1$,
every smoooth, or nodal, element in  $L_g$ is a Brill-Noether curve and the general element in $L_g$ is a Petri curve \cite{ABFS16}.

 \item Consider the blow-up of $Y$ in $p_{10}(g)$. Denote by $\Gamma_g$ the proper transform of $L_g$ under this blow-up and by $J$ the proper transform of $J_Y$.
 Then, blow-down  $J$ via $\Gamma_g$. The resulting polarized surface $(\ov X, \ov C)$, is a surface in $\PP^g$ of degree $2g-2$ with canonical sections and one elliptic singularity as its only singularity.

\item It was proved in \cite[Theorem 10.3, p.583]{ABS17} that $(\ov X, \ov C)$ is a limit of smooth polarized K3 surfaces $(S,C)$.
\end{enumerate}}

Consider the diagram
\be\label{diagramX}
\xymatrix{X\ar[r]^{\sigma\,\,\,\,\,\,\,}\ar[d]_\eta&\ov X\subset\PP^g\\
Y\ar[r]^\pi &\PP^2\\
}
\ee

where $\pi$ is the blow-up in $p_1,\dots,p_9$ and $\eta$ is the blow up of $Y$ in $p_{10}(g)$. The morphism $\sigma$ is the blow-down of $J$. 
We let  $\ell,E_1,\dots, E_9\subset X$ denote the preimages under $\eta$
of the curves with the same name in $Y$,
while $E_{10}\subset X$ is the exceptional divisor of $\eta$,  
and $C\subset X$ the proper transform of a curve $C'\in L_g$ and, as we already mentioned, we let $\Gamma_g=|C|$ denote the proper transform  
of the linear system $L_g=|C'|$. We will also denote by $\ov\Gamma_g$ the image of $\Gamma_g$ under $\sigma$.

\be\label{def-X}
\text{\parbox{.88\textwidth}{{\it From now on, by abuse of language,
we will refer also to the surface $X$   as a Halphen surface and to $\Gamma_g$ as  a du Val
linear system on $X$.  }}}
\ee

Set 
$$
J'=\eta^{-1}(J_Y)=J+E_{10}.
$$
We have:
\be\label{curvX}
\aligned
C&=3g\ell-gE_1-\cdots-gE_8-(g-1)E_9-E_{10}\,,\\
&=gJ+E_9+gE_{10}\\
&=gJ'+E_9-E_{10}\\
K_X&=-J\,,\quad J^2=-1\,,\quad J'\cdot J=0\,,\quad {J'}^2=0\\
J\cdot C&=0\,,\quad J'\cdot C=1\,, \quad C_{|J}=\cO_J\,,\quad \\
\endaligned
\ee
The {\it index of $X$} is the minimal integer $k$ such that $|kJ'|$ is  a pencil.
The morphism $\sigma: X\to \ov X\subset \PP^g$ is  given by the linear system $|C|$ and  contracts $J$ to the elliptic singularity of $\ov X$.
A Halphen surface $X$, in some respects behaves like $\PP^2$,
and in some others like a K3 surface, but over both it has the advantage of
having, for every value of $g\geq2$, a canonical polarization of genus $g$: the du Val linear system $L_g$.

\section{Clifford index of du Val curves}\label{cliff-duv}

Let $A$ be a line bundle on a smooth genus-$g$ curve $C$. 
Recall the definitions of Clifford index of $A$ and Clifford index of $C$:

$$
\cliff(A)=d(A)-2r(A)\,,\qquad\text{where}\,\, d(A)=\deg A\,,\quad r(A)=h^0(A)-1\,,
$$
so that
$$
\cliff(A)=g+1-h^0(A)-h^1(A)\,.
$$
$$
\cliff(C)=\min\{\cliff (A)\,|\, h^0(A)\geq 2\,,h^1(A)\geq 2\}
$$
The following inequality holds for every curve $C$:
$$
\cliff(C)\leq\left\lfloor{\frac{g-1}{2}}\right\rfloor
$$
and for a general curve of genus $g$, we have $\cliff(C)=\lfloor{\frac{g-1}{2}}\rfloor$.
Denoting by $\gamma(C)$ the gonality of $C$, it is well known that
$$
\gamma(C)-3\leq\cliff(C)\leq \gamma(C)-2
$$

A remarkable theorem by Green and Lazarsfeld \cite{GL87} shows that the smooth elements belonging to any given  ample linear system on a K3 surface all have the same Clifford index, denoted by $c$, and that there exists a line bundle $L$ on $X$, with $\deg L_{|C}\leq g-1$, such that, for every smooth element $C\in \Gamma_g$,
\be\label{GLlb}
\cliff(L_{|C})=c\qquad (\deg L_{|C}\leq g-1)
\ee
  The aim of this note is to show that, also from the point of view of the Clifford index, du Val systems on Halphen surfaces behave like hyperplane sections of a K3 surface, namely:

%MAIN THEOREM
 \begin{thm}\label{main} Let $\Gamma_g$ be the  du Val linear system of genus $g\geq2$ on a Halphen surface $X$. Then 
 
 i) The smooth curves $C$ in $\Gamma_g$ have all the same Clifford index  $c$.

ii) If
 $
c<\left\lfloor{\frac{g-1}{2}}\right\rfloor
$, and $k=\ind(X)$, 
then the smooth curves in $\Gamma_g$ also have the same gonality and
$$
 \gon(C)=\cliff(C)+2=k
 $$
and  the elliptic pencil $|kJ'|$ cuts out on $C$ a $g^1_k$ which is totally ramified in   $E_{10}\cap C$. 
Moreover  the Green-Lazarsfeld's line bundle $L$ in (\ref{GLlb})
is the elliptic pencil $\cO_X(kJ')$, except, perhaps, in genus $g=5$ where one could have 
$\cO_X(4J')=|C-L|$.

 \end{thm}

 In proving part i) of this theorem
 we follow a path that is completely parallel to the one indicated by Green and Lazarsfeld
but that, due the different nature of the surfaces involved (K3 surfaces on the one side,  and Halphen surfaces on the other),
 must occasionally diverge into some twists and turns  before reaching its natural  destination. 
The main difference from the K3 case is, of course, the non-vanishing of the canonical sheaf for Halphen surfaces, affecting both Serre duality and the Riemann-Roch theorem. In addition, the main properties of the Lazarsfeld-Mukai bundles, in the two cases, are somewhat different.
 The reader of the present paper should keep on hand the paper by Green and Lazarsfeld's, for consultation.
\vskip 0.2 cm
In proving part ii) of this theorem we will make essential use of the equality: $J'\cdot C=1$.

\section{Some tools}

 We will freely use Green-Lazarsfeld's Proposition (1.1), regarding smooth rational surfaces $S$ with vanishing $h^1(\cO_S)$. 
In \cite[Proposition (1.5)]{GL87} the authors show that, when dealing with a K3, 
their Proposition 1.1 holds 
 with the weaker assumption that the vector bundle $E$ is generated by its sections outside a finite set of points.
In the case of a Halphen surface $X$  \cite[Proposition (1.5)]{GL87} should be changed into the following.
  
 \begin{prop}\label{ex1.5} Let $E$ be a vector bundle on a Halphen surface $X$. Assume that $H^2(E)=H^0(E^*)=0$, and  that
$E$ is generated by its global sections away from finitely many points {\rm not belonging to $J$}.  Then:

{\rm (i)}  If $c_1(E)^2>0$, then $h^0(E)\leq h^0(\det E)$.

{\rm(ii)}  If $c_1(E)^2=0$, and $H^1(E^*)=0$, then  $E=\oplus^n\cO_X(\Sigma)$ where $\Sigma$ is a smooth irreducible curve on $X$ which moves in a base-point free pencil.
 \end{prop}

 \vskip 0.2 cm 
 
 The proof of this Proposition is based on the following variation of Lemma (1.6) in \cite{GL87}.
 
 \begin{lemma}\label{16v} Let  $E$ be a vector bundle on a Halphen  $X$ which is generated by its global sections away from finitely many points {\rm not belonging to $J$}, and such that $H^0(E^*)=0$. Then there exists a globally generated vector bundle $F$ on $X$ with $\det F=\det E$, $h^0(F)\geq h^0(E)$, and $h^2(F)=h^2(E)=0$.
  \end{lemma}
  The request that the zero-dimensional scheme, over which $E$ is not generated by its sections, should be disjoint from $J$, plays a central role  throughout  this note.
 Here, we point out the variations to be made in the proof of Lemma \ref{16v} with respect to the K3 case. 
 
 \proof (of Lemma \ref{16v}) The exact sequence   \cite[(1.7)]{GL87} stays, with the added condition that $\supp(S_E)\cap J=\emptyset$
 (since by hypothesis $E$ is generated along $J$). From this and from  the exact sequence  \cite[(1.8)]{GL87} we deduce the exact sequence
 $$
 0\longrightarrow E^*(-J)\longrightarrow H^0(E)^*\otimes\cO_X(-J)\longrightarrow V^*(-J)\longrightarrow R_E\longrightarrow 0\,.
 $$
We  observe again  a recurrent trait in this note:  tensoring by $\cO_X(-J)$ should never effect the relevant  torsion sheaves.
By Serre duality  the map $H^2(E^*(-J))\to H^0(E)^*\otimes H^2(\cO(-J))$ is an isomorphism, and this implies that $H^0(V^*(-J))\to H^0(R_E)$ is surjective, and so is, {\it a fortiori},  $H^0(V^*)\to H^0(R_E)$. As in  \cite{GL87}, the global generation of $F$ follows from this surjection.
Then  the exact sequence    \cite[(1.9)]{GL87}, shows that: "$H^i(E^*)=0\,\Rightarrow H^i(F^*)=0$", for $i=0,1$, and that $\det(F)=\det(E)$. 
The exact sequence    \cite[(1.9)]{GL87}, 
 twisted by $\cO_X(-J)$, together with Serre duality  show that   $h^0(F)\geq h^0(E)$, and $h^2(F)=h^2(E)=0$.
  \endproof

  \proof (of Proposition \ref{ex1.5}): The case (i) is treated exactly as in the proof of   \cite[Proposition (1.5)]{GL87}. In case (ii),  following again the proof of\cite[Proposition (1.5)]{GL87},
  we are reduced to the case where there is  a short exact sequence $0\to H\to \oplus\cO_X(\Sigma)\to E\to 0$, with $H$ a trivial bundle, and $\Sigma$ a smooth irreducible curve moving  in a base point free pencil.
  Twisting this sequence by $\cO_X(-J)$, and taking into account that, by  hypothesis,
  $H^1(E^*)=0$, 
  we deduce that $H=0$.\endproof
  \skip 0.3 cm
    \begin{rem}\label{ellpt-penc} Regarding the base-point-free pencil $|\Sigma|$ in item {\rm (ii)} of the preceding Proposition, notice 
  that:

  a) if $\Sigma\cdot J=0$, the genus formula gives $g(\Sigma)=1$.

  b)  if $\Sigma\cdot J'=0$, the HIT, gives  
$|\Sigma|=|kJ'|$ (see point (v) in Section 2),
  where $k$ is the index of $X$.

  \end{rem}

We now turn our attention to Lazarsfeld-Mukai bundles in the context of Halphen surfaces.
Let $A$ be a base-point-free line bundle on a smooth curve $C\in |\Gamma_g|$.
Consider the Lazarsfeld-Mukai bundle $E_{A,C}$ defined by either of the following two sequences

\be\label{LM}
0\longrightarrow E_{A,C}^*\longrightarrow H^0(C,A)\otimes\cO_X
\longrightarrow \iota_*A
\longrightarrow 0
\ee

\be\label{LM-dual}
0\longrightarrow H^0(C,A)^*\otimes\cO_X\longrightarrow E_{A,C}
\longrightarrow \iota_*(K_CA^{-1})
\longrightarrow 0
\ee

 From now on, we will omit the symbol $\iota_*$ and we will  write $E_A$ instead of $E_{A,C}$ .
Here is a list of some straightforward computations.

\be\label{comp}
\aligned
&h^i(E_A)=h^{2-i}(E_A^*(-J)))\,,\\
&h^0(E_A)=h^0(A)+h^1(A)\,,\quad H^1(E_A)=H^0(A)^*\,,\,\quad(\textcolor{red}{\dagger})\,\,\,\quad H^2(E_A)=0\,,\\
&H^0(E_A^*)=0\,,\quad H^1(E_A^*)=0\,,\quad h^2(E_A^*)=h^1(A)\,.\,\,\quad(\textcolor{red}{\dagger})
\endaligned
\ee

\be\label{c2}
c_1(E_A)=[C]\,,\qquad c_2(E_A)=d=\deg A\,.
\ee

\be\label{cJ}
C\cap J=\emptyset, \quad\quad\text{so that}\quad c_1(E_A)\cdot J=0\,.
\ee

As a consequence:
\be\label{restrJ}
\text{\it The  restriction of $E_A$ to $J$ is trivial.}
\ee

 We also have:

\be\label{chi-ea}
\chi(E_A\otimes E_A^*)=-h^0(A)^2+2(1-\rho(A))\,\,\quad(\textcolor{red}{\dagger})\footnote{ 
\quad When $X$ is a K3 surface, in the three cases marked with  a dagger $(\textcolor{red}{\dagger})$ we have, respectively, $H^1(E_A)=0$, $h^2(E^*_A)=h^1(A)+h^0(A)$, \quad$\chi(E_A\otimes E_A^*)=2(1-\rho(A))$.
}
\ee

where $\rho(A)=g-h^0(A)h^1(A)$ is the Brill-Noether number of $A$.

\section{Green and Lazarsfeld's idea}
The procedure invented by Green and Lazarsfeld to prove the constancy of the Clifford index 
on smooth hyperplane section of a given polarized K3 surface can be imitated in our setting as follows.
Start with a smooth curve $C\in \Gamma_g$ having minimal Clifford index, and let $A\in \Pic(C)$ be a line bundle realizing the Clifford index of $C$:
\be\label{Acomp}
\cliff(A)=\cliff(C)\,.
\ee
One may always assume that

\be\label{assume}
\aligned
&\deg A\leq g-1\,,\\
&\cliff(C)<\Bigl\lfloor\frac{g-1}{2}\Bigr\rfloor\,,\\
&\ind(X)-2\geq\cliff (C)\,,\,\, \text{(otherwise $\cliff(C)>\gamma(C)-2$)}\,,\\
&\text{both $A$ and $K_CA^{-1}$, and therefore $E_A$, are generated by their global sections.}
\endaligned
\ee

To prove part i) of Theorem \ref{main}, it suffices to produce a line bundle $L$ on $X$
that contributes to the Clifford index of every smooth curve $C_0\in  \Gamma_g$ (i.e. $h^0(L_{|C_0})\geq 2$, $h^1(L_{|C_0})\geq2$)
and such that 
\be\label{e}
\cliff (L_{|C_0})\leq\cliff(C).
\ee
Look at the Lazarsfeld Mukai bundle $E_A$. Ideally one would like to find a line bundle $N$, with $h^0(N)\geq2$,  and an exact sequence of sheaves
\be\label{seq-N}
0\longrightarrow N\longrightarrow E_A\longrightarrow F\longrightarrow \tau\longrightarrow0
\ee
where $F$ is locally free, and $\tau$ supported on a finite set of points {\it not in $J$}. Let's see why this would be ideal.
The basic observation regarding a sequence like (\ref{seq-N}) is the general Lemma \cite[Lemma 3.1]{GL87}
which we write down,  for the convenience of the reader, with some  necessary modifications and additions with respect to the original one.

\begin{lemma} \label{31}Let $E$ be a rank-$n$ vector bundle on a Halphen surface $X$, which is generated by its global section away 
from a finite set of points {\rm not including points in $J$}. Suppose we are given an exact sequence 
\be\label{seq-N2}
0\longrightarrow N\longrightarrow E\longrightarrow F\longrightarrow \tau\longrightarrow0
\ee
where $N$ is a line bundle with $h^0(N)\geq2$, $F$ is locally free, and $\tau$ is supported on a finite set {\rm not including points in $J$}.
Assume that $H^2(E)=H^0(E^*)=H^1(E^*)=0$.  Let $\Gamma$ be a divisor representing $c_1(E)$.
Then

{\rm (i)} $\dim |\Gamma|\geq 1$, and $|\Gamma|$ has no fixed components.

{\rm (ii)} $F$ is generated by its global sections off a finite set not including points of $J$, and $\cO_X(\Gamma)=N\otimes\det F$.
Furthermore: 
\be\label{vanish}
\aligned
&h^2(F)=h^0(F^*)=h^1(F^*)=0\,,\\
&h^1(F)\leq h^1(E)\,,\\
&h^0(\det F)\geq 2.
\endaligned
\ee
{\rm (iii)} If $C\subset X$ is a curve such that $C\cap J=\emptyset$, and $C-\Gamma$ is effective, then
$$
h^0(N_{|C})\geq h^0(N)\geq 2\,,\quad h^1(N_{|C})\geq h^0(\det F)\geq 2\,.
$$
\end{lemma}
 The proof of this Lemma, which  is left to the reader, needs only a few modifications from the original one. For example, on p. 368, l. +9 
the sentence ``{\it but $u$ cannot be an isomorphism since $H^2(E)=0$...}" should be substituted by: 
"{\it but $u$, and therefore $u(-J)$, cannot be an isomorphism since $H^2(E(-J))=0$...}". Also $h^2(N)=h^0(N^*(-J))=0$.
\vskip 0.2 cm 
Let us go back to the ``ideal situation" where, for some reasons, we find an exact sequence like  (\ref{seq-N}) directly involving the Lazarsfeld -Mukai bundle $E_A$, and let us show that, at least when $c_1^2(F)>0$,  we immediately get a line bundle $L$ satisfying (\ref{e}), thus proving Theorem \ref{main}.
In fact,  we have
\be\label{c1^2}
\aligned
g+1-\cliff(C)&=g+1-\cliff(A)\\
&= h^0(E_A)\\
&\leq h^0(N)+h^0(F)\\
&\leq h^0(N)+h^0(\det F) \quad \text{[by i) of Proposition \ref{ex1.5}]}\\
&\leq h^0(N_{|C})+h^1(N_{|C})\quad  \text{[by (iii) of Lemma \ref{31}]}\\
&\leq g+1-\cliff(N_{|C})\\
\endaligned
\ee
and we are done by setting $N=L$.  The case $c_1(F)^2=0$ is more delicate and requires  a subtle use of part ii) 
of Proposition \ref{ex1.5}. Regardless of this added difficulty, things  are anyhow more complicated and the ideal situation does not necessarily occur.
Green and   Lazarfeld  invented a reduction procedure to "cut down"  a Lazarsfeld-Mukai bundle, without losing some key
information, to get  lower rank vector   bundles $E$ to which to apply Lemma \ref{31}. In order to be able to initiate the above string  
of inequalities these bundles should certainly satisfy the condition: $h^0(E)\geq h^0(E_A)$.

\section{The reduction procedure}

The following definition of {\it reduction} is a slight variation of   \cite[Definition 2.7]{GL87} that is  suited to Halphen surfaces
[see conditions R0) and R1) below].
 
\begin{defin}\label{red}
Let $E_0$ be a vector bundle on a Halphen surface $X$. 
A vector bundle $E_1$ of rank 
$\geq 2$ is a reduction of $E_0$, if the following properties are satisfied.

{\rm R0)} $E_1$ satisfies the  following three conditions:   
\hskip 0.8 cm 

\hskip 0.3 cm{\rm a)} 
 $h^0(E_1^*)=h^1(E_1^*)=h^2(E_1)=0$.

\hskip 0.3 cm{\rm b)} The restriction of $E_1$ to $J$ is trivial.

\hskip 0.3 cm{\rm c)} $c_1(E_1)$ is represented by an effective divisor.
\hskip 0.3 cm

{\rm R1)} There is a map $E_0\to E_1$ which is surjective off a finite set $T$, with $T\cap J=\emptyset$.

{\rm R2)} $h^0(E_1)\geq h^0(E_0)$.
 
{\rm R3)} $\det E_1=\det E_0\otimes \cO_X(-D)$, for some effective or zero divisor $D$.

{\rm R4)} $c_1^2(E_1)-4c_2(E_1)+8r(E_1)\geq c_1^2(E_0)-4c_2(E_0)+8r(E_0)$.
\end{defin}
 Observe that R0) b) implies that
 \be\label{dot-J}
 c_1(E_1)\cdot J=0
 \ee
Conditions a), b), c) in R0) are satisfied by any  Lazarsfeld-Mukai bundle on a
Halphen surface [cf. (\ref{comp}), (\ref{c2}), (\ref{restrJ})], and among them
 property R0) b) is crucial.
Also observe that, since $C\cap J=\emptyset$, a Lazarsfeld-Mukai vector bundle $E_A$, as above, is a reduction of itself
and that in fact any vector bundle satisfying R0) is a reduction of itself. A reduction is {\it minimal} if there is no reduction of strictly lower rank. Clearly reductions of reductions are reductions, and  minimal reductions exist  even if, in our treatment, we will not make use of this last fact. By R3), for any reduction $E$ of a
 Lazarsfeld-Mukai bundle $E_A$, the class $c_1(\det E)$ is represented by a sub-divisor of a divisor in $L_g=|C|$.

\vskip 0.2 cm
Next, notice that, for a vector bundle $F$ on $X$, we have:
$$
\chi(F\otimes F^*)=r(F)^2+(r(F)-1)c_1(F^2)-2r(F)c_2(F).
$$
We also set
$$
(R4)(F)=c_1^2(F)-4c_2(F)+8r(F)
$$
Thus, when $r(F)=2$, we get

\be\label{chif}
(R4)(F)=\chi(F\otimes F^*)+12
\ee
The reduction procedure introduced by Green and Lazarsfeld at each stage  makes
the rank  drop by one and  the second Chern class increase by at least two (see (\ref{rank-dec})), so that the rationale underlying condition R4) is  to keep
all the invariants under control.

\vskip 0.2 cm 
The statement of   \cite[Proposition 2.8]{GL87} is slightly changed but the gist is the same, namely:

\begin{prop}\label{ex-prop-28}
Start with a Lazarsfeld-Mukai bundle $E_A$ as above [see (\ref{assume}]. 
Then there exists a reduction $E$ of $E_A$ equipped
an exact sequence (\ref{seq-N2}) as in Lemma \ref{31}. Furthermore, since $E$ is a reduction,  we have $h^2(E)=h^0(E^*)=h^1(E^*)=0$, and 
therefore conditions {\rm (i), (ii), (iii)} of Lemma \ref{31} are satisfied.
\end{prop}

\proof We start with $E_A=E_0$ and we will recursively produce vector bundles $E_i$, $i=1,\dots, n$, where $E_i$ is a reduction of $E_{i-1}$ and 
\be\label{rank-dec}
\rk E_i=\rk E_{i-1}-1\,,\quad i=1,\dots, n
\ee
and such that $E=E_n$ satisfies the requirements of  Lemma \ref{31}. We will follow the ideas contained in the 
   proof of  \cite[Proposition 2.8] {GL87}.   Let us first assume that there 
is a reduction $E$ of $E_A$ with 
$\rk(E)=2$, a condition that, {\it per se}, makes this reduction automatically minimal, and let us prove that $E$ satisfies the requirements of  Lemma \ref{31}.
The first task is to  prove that
\be\label{end}
h^0(E\otimes E^*)\geq 2
\ee
Suppose not, then
 condition R4), and (\ref{chif}) give:
$$
(R4)(E)=\chi(E\otimes E^*)+12=1-h^1(E\otimes E^*)+12\geq (R4)(E_A)=2g-2-4\cliff(C)+8
$$
As $\cliff(C)<\lfloor{\frac{g-1}{2}}\rfloor$, we have
$1-h^1(E\otimes E^*)\geq0$.
On the other hand, always under the assumption that $h^0(F\otimes F^*)=1$, the standard exact sequence describing the restriction of $E\otimes E^*$ to $J$, and the fact that, by hypothesis, this restriction
is trivial,  give: $h^1(E\otimes E^*)\geq2\rk(E)-1=3$, which is absurd. Therefore there is a non-trivial endomorphism $\phi:E\to E$, and in order to find the line bundle $N$ and the exact sequence (\ref{seq-N2}), one proceeds exactly as in \cite[p. 365, last paragraph]{GL87}, with the added remark that, being $E$ trivial along $J$, we may say that the map $E\to N$, constructed in that paragraph, is surjective away from finitely many points none of which belong to $J$.
\vskip 0.2 cm 
To construct the sequence of reductions $E_1,\dots,E_n$ of $E_A$ we  proceed inductively. 
Let $E=E_i$ be a reduction of $E_A$ and assume, as we now may, that its rank is
$\geq 3$. We want to construct $F=E_{i+1}$.
The first step  is 
 to find two points $x$ and $x'$ in $X$ such that $ H^0(E\otimes I_{\{x,x'\}})\neq 0$.  
 This is done
exactly as in  \cite[p. 366, point  2.10,  first three paragraphs]{GL87}, where the formula
$$
h^0(E)\geq 2\rk(E)
$$ 
is highlighted. Since  in  the reduction process we are about to describe the rank drops by one at each stage we have
\be\label{alternative}
\text{either}\,\,\, E=E_A\,, \,\text{or}\,\,\, h^0(E)>2rk(E). 
\ee
In  case  $h^0(E)>2\rk(E)$ one can clearly choose $x$ and $x'$ to be generic.  
Suppose $h^0(E)=2\rk(E)$, so that $E=E_0=E_A$.
Following  \cite[p. 366 ]{GL87}, one can choose $x$ generic, then
look at the natural map $\nu_x: H^0(E_A\otimes I_x)\otimes\cO_X\to E_A$ and take $x'\in \Gamma_x=(\det\nu_x)_0\in |\det E_A|=|C|$, and $s\in\nu_x^{-1}(x')$. Since  no element $ \Gamma_x\in |C|$ can be a multiple of $J$,
we may choose $x,x'\notin J$ (in fact). But we have to do more.

\begin{lemma}We may choose $x$ and $x'$
in such a way that there is a section $s\in H^0(E\otimes I_{\{x,x'\}})$ whose zero set  is disjoint from $J$. 
\end{lemma}
\proof Suppose we can't, this means
that, for all pairs of points $x, x'\in X$
\be\label{incl-xx'}
H^0(E\otimes I_{\{x,x'\}})\,\,\subset \,\,\underset{y\in J}\cup H^0(E\otimes I_{y})\
\ee
Since $E$ is trivial along $J$, and  $h^1(E^*)=0$, one sees that, for $y\in J$, we have $h^0(E\otimes I_y)=h^0(E(-J))=h^0(E)-\rk(E)$, so that
$H^0(E\otimes I_y)=H^0(E(-J))$. Thus relation (\ref{incl-xx'})  becomes: $H^0(E\otimes I_{\{x,x'\}})\subset H^0(E(-J))$ which implies:
$$
0\neq H^0(E\otimes I_{\{x,x'\}})\subset H^0(E(-J)\otimes I_{x})\,.
$$
But, since  $h^0(E)\geq2\rk E$, for general $x$ we have 
$h^0(E(-J)\otimes I_{x})=h^0(E)-2\rk E$, which is a contradiction if $h^0(E)=2\rk E$. When $h^0(E)>2\rk E$, then we can take both $x$ and $x'$ general and we would get the equality $H^0(E\otimes I_{\{x,x'\}})=H^0(E(-J)\otimes I_{x})$, for $x'$ general,  but 
then both terms of the last equality would be equal to $0$,
again a contradiction.
\endproof

Using the section $s$ constructed  in the Lemma, we proceed as in \cite[p. 366] {GL87} to find a divisor $\Delta$, and   a zero-dimensional subscheme
$Z\subset X$ such that $(s)_0=Z\cup \Delta$, with $\supp(Z\cup\Delta)\cap\supp J=\emptyset$, and get dual exact sequences 
\be\label{ex-seq-N*}
0\longrightarrow F^*\longrightarrow E^*\longrightarrow\cO_X(-\Delta)\longrightarrow\cO_Z(-\Delta)\longrightarrow0
\ee

\be\label{ex-seq-N}
0\longrightarrow\cO_X(\Delta) \longrightarrow E\longrightarrow F \longrightarrow\tau \longrightarrow0\,,\qquad 
\ee
defining the locally free sheaf $F$.  Notice that
\be\label{nt-supp} 
\supp\tau
=\supp Z\,,\qquad \supp\tau\cap  \supp J=\emptyset
\ee

Three possibilities may arise:

\item i) $\Delta=0$,

\item  ii) $\Delta\neq0\,,\text{but}\,h^0(\cO_X(\Delta)=0$,

\item iii) $h^0(\cO_X(\Delta)\geq 2$.

In case (iii) we  set $N=\cO_X(\Delta)$ and we get an exact sequence (\ref{seq-N2}) as in Lemma \ref{31}, so we are done.
In  contrast with the K3 case, let us show that
 ii) never occurs: indeed,  assume that $\Delta$ is irreducible. Since $\Delta\cdot J=0$, from the Riemann-Roch theorem we get 
 $-2h^1(\Delta)=\Delta^2$. Since $X$ is unnodal we must have 
 $h^1(\Delta)=0$. On the other hand, tensoring with $\cO_X(-\Delta)$ the sequence $0\to\cO_X(-J)\to\cO_X\to\cO_J\to 0$,
we get an injection $H^0(\cO_J)\to H^1(-\Delta-J)= H^1(\Delta)^*$, which is absurd.
\vskip 0.2 cm 
In case i) we claim that $F$ is a reduction of $E$. In this case
(\ref{ex-seq-N*}) and (\ref{ex-seq-N}) read:
\be\label{ex-seq-N*2}
0\longrightarrow F^*\longrightarrow E^*\longrightarrow\cO_X\longrightarrow\cO_Z\longrightarrow0
\ee

\be\label{ex-seq-N2}
0\longrightarrow\cO_X \longrightarrow E\longrightarrow F \longrightarrow\tau \longrightarrow0\,,\qquad 
\ee
 We must verify R0)-R5) , where the pair $(E_0,E_1)$ in Definition \ref{red} is now the pair $(E,F)$.
 Properties R1) and R3) follow from 
(\ref{ex-seq-N2})
 and the fact that $\supp \tau\, \cap\,\supp J=\emptyset$. Next, since $Z$ contains at least $x$ and $x'$, one easily computes: $h^0(\I_Z)=h^2(\I_Z)=0$,
 $h^1(\I_Z)\geq1$, so that $h^0(F^*)=h^1(F^*)=0$. Also from (\ref{ex-seq-N2}) we get $h^2(F)=0$, and $\det(F)=\det E$.  Thus R0) a)
 and R0) c) are satisfied. Also, since $\supp \tau \cap\supp J=\emptyset$, we get the exact  sequence 
 $0\to\cO_J\to E_J\to F_J\to 0$ and since $ E_J$ is trivial, $ F_J$ is also trivial, proving R0) b). To prove R2), using again the fact  
 that $\supp \tau \cap\supp J=\emptyset$,  we get the exact sequence $0\to F^*(-J)\to E^*(-J)\to I_Z(-J)\to 0$
 and the vanishing of $h^2(I_Z(-J))=0$. Thus the surjection $H^2(F^*(-J))\to H^2(E^*(-J))\to 0$, gives by duality $h^0(F)\geq h^0(E)$, proving  R2). We can then start again with the reduction $F$ and repeat the above procedure to reach either a reduction of rank 2
 or one for which iii) holds.

 \endproof
 
 \begin{rem}\label{steps}  $\phantom{x}$
 
 {\bf a)} Notice that, by construction,   $\Delta$ is an effective subdivisor of $c_1(E)$  and, since $E$ is a reduction of $E_A$,
 by {\rm R0) c)} and {\rm R3)}  also $[C-\Delta]$ is effective.
 \vskip 0.2 cm 
 {\bf b)}
 From the proof of Proposition \ref{ex-prop-28} we see how the reduction process works.
 We start with $E_A:=E_0$, and we construct vector bundles $E_1,\dots, E_n$, each a reduction of the preceding one, with ranks decreasing by one at each step,
 and equipped with exact sequences
 \be\label{passi}
 \xymatrix{
 0\ar[r]&\cO_X\ar[r]& E_0\ar[r]^{\alpha_1}&E_1\ar[r]&\tau_1\ar[r]  &0\\
 &&\cdots\cdots\cdots\\
  0\ar[r]&\cO_X\ar[r]& E_{n-2}\ar[r]^{\alpha_{n-1}}&E_{n-1}\ar[r]&\tau_{n-1}\ar[r]  &0\\
   0\ar[r]&\cO_X(\Delta)\ar[r]\ar@{=}[d]& E_{n-1}\ar[r]^{\alpha_n}\ar@{=}[d]&E_n\ar[r]\ar@{=}[d]&\tau_n\ar[r]\ar@{=}[d]  &0\\
 0\ar[r]&N\ar[r]& E\ar[r]&F\ar[r]&\tau\ar[r]  &0\\
 }
 \ee

The $\tau_i$'s are supported on points of $X$  {\rm away from $J$}, and also
 \be\label{int-J}
E_i\cdot J=0\,, \quad\text{for}\,\,\,\,i=0,\dots, n-1\,,\quad
 \Delta\cdot J=0\,,\,\, \text{so that}\,,\quad c_1(F)\cdot J=0
 \ee
Finally $h^0(N)\geq2$, and since $N\cdot J=0$, the linear system $|N|$ has no base point on $J$, i.e. $N$ is generated by its sections along $J$. Finally, by {\bf a)}, both $N$ and $C-N$ are effective.

 \end{rem}

 We now come to the proof of the main Theorem  \ref{main}.
 
 \proof (of  {\bf part (i) of Theorem \ref{main}}). The proof runs exactly as the proof of Green and Lazarsfeld's Theorem \cite[pp. 368-370]{GL87}. One starts with a smooth curve $C\in \Gamma_g$  and a line bundle $A\in \Pic(C)$ satisfiying (\ref{Acomp}) and (\ref{assume}), and  looks at the Lazarsfeld-Mukai bundle $E_A$. By Proposition \ref{ex-prop-28} there exists a reduction $E$ of $E_A$ equipped with an exact sequence  (\ref{seq-N2}) enjoying properties (i),(ii), and (iii) of Lemma \ref{31}.

{\bf Case 1:} $c_1(F)^2>0$. This is dealt with as in (\ref{c1^2}) \cite[p. 369]{GL87};

{\bf Case 2:}  $c_1(F)^2=0$. This case is more delicate.  By virtue of (\ref{vanish}), we can apply Proposition \ref{ex1.5} to the vector bundle $F$ to obtain that $F=\oplus \cO_X(\Sigma)$ where $|\Sigma|$ is a base point free pencil on $X$. By the last observation in Remark \ref{steps}, we have that $N\cdot J=0$; it then follows that $\Sigma\cdot J=0$, as well. Thus, by Remark \ref{ellpt-penc},
the pencil $|\Sigma|$ is an elliptic pencil.
Next, exactly as in  \cite[p. 369]{GL87},
one shows 
 that there is an exact sequence
 
\be\label{ex3.5}
0\longrightarrow N_0\longrightarrow E_{KA^{-1}}\longrightarrow F_0\longrightarrow\tau_0\longrightarrow0
\ee
satisfying the hypotheses of Lemma \ref{31}.
However, in our context, it is necessary to prove, in addition,  that {\it none of  the points in support of $\tau_0$  lies in $J$}.
For this  we  recall some of the details of  the construction of the sequence (\ref{ex3.5}).
One considers the diagram  \cite[p. 369]{GL87}

\be\label{tau-zero}
\xymatrix{0\ar[r]&E^*_{KA^{-1}}\ar[r]\ar[d]^u&H^0(E_A)\otimes\cO_X\ar[r]\ar[d]^w&E_A\ar[r]\ar[d]^\mu&0\\
0\ar[r]&\cO_X(-\Sigma)\ar[r]&H^0(\Sigma)\otimes\cO_X\ar[r]&\cO_X(\Sigma)\ar[r]&0\\
}
\ee
where the first row is Tjurin's exact sequence, and $\mu$ is the composition of $\alpha_1\circ\cdots\circ\alpha_n: E_A\to F$ [see (\ref{passi})] with the projection onto one of the summands of $F=\oplus\cO_X(\Sigma)$. If $D$ is the largest effective or zero divisor such that $u$ factors through $O_X(-\Sigma-D)$, 
one gets the exact sequence defining $F_0$:
\be\label{ex3.5}
0\longrightarrow F_0^*\longrightarrow E_{KA^{-1}}^*\overset v\longrightarrow O_X(-\Sigma-D)\longrightarrow O_Z(-\Sigma-D)\longrightarrow0
\ee
where $Z$ is a zero-dimensional scheme with the same support as $\tau_0$.  The proof that the support of $\tau_0$ is disjoint from $J$,
follows by chasing  diagram (\ref{tau-zero}), using: the snake lemma, the fact that the support of the $\tau_i$'s, in diagram (\ref{passi}), are disjoint from $J$, and the fact that $N$ (in diagram (\ref{passi})) is generated by its global sections outside the support of $J$ . The rest of the proof of part (i) of Theorem \ref{main},  is identical to    \cite[end of p. 369, and p. 370]{GL87}.  We recap the resulting alternatives.
\vskip 0.2 cm
Start with $0\to N\to E\to F\to \tau\to 0$. A priori the following cases are possible:

{\bf a)} If $c_1(F)^2>0$, take $L=N$, and therefore $M=c_1(F)$. 

{\bf b)} If $c_1(F)^2=0$, there exists   an exact sequence $0\to N_0\to E_{K_CA^{-1}}\to F_0\to \tau_0\to 0$ and:

{\bf b1)} if $c_1^2(F_0)>0$, we take $L=c_1(F_0)$, $M=N_0$,

{\bf b2)} if $c_1^2(F_0)=0$, there exists an elliptic pencil $\Sigma_0$ such that $F_0=\cO(\Sigma_0)$, then we take $L=F_0$ and then $M=N_0$
\vskip 0.2 cm

We now prove {\bf  part (ii) of Theorem \ref{main}}.  Here it is interesting to compare our arguments
to  the ones in  \cite{Mart89}, and in \cite{Knut01}. 
By part (i) all smooth element $C\in \Gamma_g$ have the same Clifford index $c$, and there is 
 a line bundle $L$ on $X$ such that, for all smooth elements $C\in \Gamma_g$
$$
\cliff(C)=\cliff(L_{|C}):=c\,,\quad 0\leq c<\left\lfloor{\frac{g-1}{2}}\right\rfloor
$$
so that in particular $h^0(L_{|C})\geq 2$, and $h^1(L_{|C})\geq 2$, and we set as usual
$M=C-L$. We claim that
\be\label{complete}
h^0(L)=h^0(L_{|C}).
\ee
This is true by a straightforward adaptation of Martens' argument   \cite[proof of Lemma 2.2 ]{Mart89},   
unless we are in the case ${\bf b2)}$ with $c=0$. Then
we can take $L$ as an elliptic pencil $\cO(\Sigma_0)$.
Since $c=0$, we could have started with the pencil $|A|=g^1_2$ and 
(\ref{complete}) is obvious.

\vskip 0.2 cm
Since $J$ is the only positive divisor contracted by $L_g$, and because of (\ref{complete}), the only base  curves that $|L|$ might have are the multiples of $J$. By subtracting them
from $L$ we may assume that $|L|$ has no base points.  Now we must consider two cases.
\vskip 0.2 cm
1) $L^2>0$
\vskip 0.2 cm
2) $L^2=0$
\vskip 0.2 cm
Suppose we are  in case 1). In this case $L$ is big and nef and therefore $H^	1(L-J)=H^1(-L)=0$, so that $h^0(M)=h^0(M_{|C})$. We claim that $M\cdot J'=0$.
We know that $(L+M)\cdot J'=1$. Suppose $M\cdot J'\neq0$. Then $L\cdot J'=0$ so that, by the Hodge index theorem, $L=\alpha J'$
against the hypothesis that $L^2>0$. Thus $M=\alpha J'$, for some positive $\alpha$. Since $h^0(M)=h^0(M_{|C})$, and $|M_{|C}|$ is base point free we must have $\alpha=mk$. But then $m+1=h^0(M)=h^0(M_{|C})$. Thus $c=mk-2m$. On the other hand,
the totally ramified $g^1_k$ cut out by $|kJ'|$ on $C$ gives: $c=k-2$, so that $m=1$, when $k\neq 2$.
When $k=2$, the elliptic pencil $|kJ'|$, cuts out  the $g^1_2$ on $C$.  Since  $\deg M_{|C}=k\geq g-1$ this case can occur only for $g=5$.

\vskip 0.2 cm
Suppose we are  in case 2).  We consider the two possibilities: either $L\cdot J'=0$, or $L\cdot J'=1$. Suppose that $L\cdot J'=1$. Then by the Hodge index theorem we have $|M|=|\alpha J'|$. But then $L=\cO(\Sigma_0)$, and since $L$ has no base curve, $h^0(L-J)=1$,
so that, by Riemann-Roch $h^1(L-J)=h^1(-L)=1$, but then $h^0(M)=h^0(M_{|C})$, and we proceed as in case 1),  showing that   $g=5$. 
We are left with the case $L\cdot J'=0$. Then $L=\alpha J'$ and, using (\ref{complete}) we proceed exactly as we did for $M$, 
proving  that $|L|=|kJ'|$.
\endproof

\vskip 0.2 cm
One final remark and a small {\it Erratum} regarding \cite {AB17}.

\begin{rem} {\rm Dealing with du Val linear systems automatically excludes the exception to the constancy of gonality in the case of K3-sections that was discovered by Donagi and Morrison \cite{DM89}, and thoroughly studied by Ciliberto and Pareschi \cite{CP95}. This exclusion of course follows from Theorem \ref{main}, but it can be easily seen directly. The du Val linear system $\ov\Gamma_2$ exhibits a Halphen surface $\ov S$ (recall diagram (\ref{diagramX})) as a double cover 
$\pi: \ov S\to \PP^2$ ramified along a sextic. The counter-example by Donagi-Morrison  looks  at the genus 10 linear system $|\pi^*(\cO_{\PP^2}(3)|$, but this is not the du Val system $\ov\Gamma_{10}$.}

\end{rem}

\vskip 0.2 cm
{\bf Filling a gap in}  \cite {AB17}. At the beginning of the proof of Theorem 4.2 of the cited paper (end of p. 283), Theorem 3 of \cite{ABFS16} is invoked. An essential hypothesis in  that theorem is that $\cliff(C)\geq3$. However,  this particular hypothesis is not checked  in the proof of the above cited Theorem 4.2. On the other hand, one of the
 hypotheses of that  Theorem is that the Halphen surface in question is of index $s+1$, where $s\geq6$. Thus Theorem \ref{main} of the present note shows that, under this hypothesis, $\cliff{C}\geq 5$, filling the gap.

%
 %%%%%%%%%%%%%%%%%%%%%%%%%%%%%%%%%%%%%%%%%%%%%%%%%

 \bibliographystyle{amsalpha}
 \bibliography{Cliff}

\begin{thebibliography}{Cliff}
 
 
 %
\bibitem[AB-17]{AB17}
Arbarello E., Bruno A.: {Rank two vector bundles on polarized Halphen surfaces and  the Gauss-Wahl map for du Val curves}. \emph{ Journal de l'\'Ecole polytechnique, Math\'ematiques}, {\bf 4}, 257--285
(2017)


 
 
 
                                                                     
 \bibitem[ABFS-16]{ABFS16}
Arbarello E., Bruno A., Farkas G., Sacc\`a G.: {Explicit Brill-Noether-Petri general curves},  
\emph{Commentarii Mathematici Helvetici} 91 (3) (2016), 477-491.

%
\bibitem[ABS-17]{ABS17}
Arbarello E., Bruno A., Sernesi E.: On hyperplane sections of K3 surfaces.  \emph{Algebraic Geometry} 4 (2017), 562-596.

\bibitem[B-LC23]{B-LC23} Bruno A., Lelli-Chiesa, M. :{Irreducibility of Severi varieties on K3 surfaces}
 \emph{arXiv:2112.09398} , Math. AG.
 
\bibitem[C-D12]{C-D12}  Cantat S.,Dolgachev I.: {Rational surfaces with a large group of automorphisms}.
\emph{Journal of the American Mathematical Society}, {\bf 25}, No. 3, 863--905 (2012).



\bibitem[CP95]{CP95} Ciliberto C.,  Pareschi G.: {Pencils of minimal degree on curves on a K3 surface}.
\emph{J, reine angew. Math.}, {\bf 460}, 15--36 (1995).


\bibitem[DM89]{DM89} Donagi R., Morison D.: {Linear systems on K3 sections}.
\emph{J. Diff. Geom.}, {\bf 29},  49--64 (1989).


 		

 \bibitem[GL87]{GL87} Green M.,  Lazarsfeld R.: {Special divisors on a K3 surface}.
\emph{Invent. math.}, {\bf 89},  357--370 (1987) . 	

\bibitem[Harb85]{Harb85}
Harbourne R.: {Complete linear systems on rational surfaces}.
 \emph{ Trans. Amer. Math. Soc.}, {\bf 289} (1), 213--226, (1985).


\bibitem[Knut01]{Knut01}
Knutsen L.A.: {On kth-order embeddings of K3 surfaces and Enriques surfaces}.
 \emph{ manuscripta math.} {\bf104}, 211-- 237,  (2001).


\bibitem[Mart89]{Mart89}
Martens G.: {On curves on K3 surfaces}.
 \emph{ Springer LNM }, {\bf 1398}, pp. 174-182
(1989).



	

\bibitem[Nag60]{Nag60}
Nagata M.:
{On rational surfaces {II}.}
\emph{ Memoirs of the College of Science, University of Kyoto},
  {\bf32}, 271--293, (1960).

\bibitem[SD74]{SD74}
Saint-Donat B.:
{Projective models of K3 surfaces}
\emph{ American Journal of Mathematics},
  { 96}, {\bf 4},  602--639, (1974).


 
 \end{thebibliography}

\end{document}